\documentclass{amsart}
\usepackage{amssymb}
\usepackage{amsfonts}
\usepackage{amstext}
\usepackage{algorithmic}
\usepackage{algorithm}
\usepackage{graphicx}
\usepackage{epstopdf}
\usepackage[all]{xy}

\allowdisplaybreaks

\numberwithin{equation}{section}
\newtheorem{theorem}{Theorem}[section]
\newtheorem{lemma}[theorem]{Lemma}
\newtheorem{proposition}[theorem]{Proposition}
\newtheorem{corollary}[theorem]{Corollary}

\theoremstyle{definition}
\newtheorem{definition}[theorem]{Definition}
\newtheorem{example}[theorem]{Example}

\theoremstyle{remark}
\newtheorem{remark}[theorem]{\bf{Remark}}

\newcommand{\Gcal}{{\mathcal{G}}}
\newcommand{\Ical}{{\mathcal{I}}}

\renewcommand{\ker}{{\rm{ker}}}

\newcommand{\tens}{\otimes}
\newcommand{\id}{{\rm id}}
\newcommand{\bn}{{}^{(0)}}
\newcommand{\bo}{{}^{(1)}}
\newcommand{\bt}{{}^{(2)}}

\renewcommand{\o}{{}_{(1)}}
\renewcommand{\t}{{}_{(2)}}

\newcommand{\eps}{\varepsilon}
\newcommand{\sig}{\sigma}

\newcommand{\la}{\triangleright}
\newcommand{\ra}{\triangleleft}
\newcommand{\iso}{\cong}
\newcommand{\lb}{\langle}
\newcommand{\rb}{\rangle}

\newcommand{\rharp}{\leftharpoondown}
\newcommand{\lharp}{\rightharpoonup}

\begin{document}

\title{Integral theory for Hopf (co)quasigroups}
\keywords{}
\subjclass{}

\author{Jennifer Klim}
\address{Queen Mary, University of London\\
School of Mathematics, Mile End Rd, London E1 4NS, UK}
\date{Version 1: \today}

\begin{abstract}
We recall the notion of a Hopf (co)quasigroup defined in \cite{Kl09} and define integration and Fourier Transforms on these objects analogous to those in the theory of Hopf algebras. Using the general Hopf module theory for Hopf (co)quasigroups from \cite{Br09} we show that a finite dimensional Hopf (co)quasigroup has a unique integration up to scale and an invertible antipode.  We also supply the inverse Fourier transformation and show that it maps the convolution product on $H$ to the product in its dual $H^*$. Finally, we further develop the theory to consider Frobenius Hopf (co)quasigroups, separability and semisimplicity.
\end{abstract}
\maketitle

\section{Introduction}

In \cite{Kl09} we defined the notion of a Hopf quasigroup as a not necessarily associative, but unital algebra with a coassociative coalgebra structure and an antipode satisfying certain conditions. We showed that a theory similar to that of Hopf algebras was possible in this case, including an application of the Hopf module lemma to study differential calculus on a Hopf coquasigroup. In \cite{Br09} the notion of Hopf modules and the Hopf module lemma or Galois property were formalized in the general case and used to characterize the axioms of a Hopf (co)quasigroup.  Here we will apply the same ideas to the notion of integration.

In Section 2 we recall the necessary background. In section 3 we define left and right integrals on Hopf quasigroups and prove some useful identities. We use the theory developed in \cite{Br09} to prove existence and uniqueness of these integrals. These definitions and properties also hold for Hopf (co)quasigroups, and this is proved in section 4. We prove existence of integrals on commutative flexible Hopf coquasigroups, by proving that if $A$ is a Hopf coquasigroup and $M$ is a right $A$-Hopf module as defined in \cite{Br09}, then $M$ is isomorphic as a right $A$-Hopf module to $M^{\widehat{coA}}\tens A$, where $M^{\widehat{coA}}$ is the set of $A$-coinvariants with respect to some induced coaction.

This leads us to define the Fourier transform $F:H \to H^*$ from a finite-dimensional Hopf quasigroup to its dual in Section 5. If we assume the Hopf quasigroup is cocommutative and flexible we obtain more familiar identities. As one would expect, the Fourier transform maps the convolution product on $H$ to the product on $H^*$. In Section 6 the Fourier transform is defined on Hopf coquasigroups in a similar manner and we consider the example of $\Gcal$ a finite (IP) quasigroup, so that, as shown in \cite{Kl09}, the group function algebra $k[\Gcal]$ is a Hopf coquasigroup. We describe the resulting integrals and Fourier transform.

In Section 7, we develop the theory further by discussing Frobenius and separable Hopf (co)quasigroups. Finally, in Section 8 we consider semisimple Hopf coquasigroups.

\section{Hopf quasigroups and Hopf modules}

We recall that an (inverse property) {\em quasigroup} or {\em IP loop} \cite{Pf90} is a set $\Gcal$ with a product, denoted by omission, an identity $e$, and for each $u\in \Gcal$ an element $u^{-1}\in \Gcal$ such that 
\[ u^{-1}(uv)=v,\quad (vu)u^{-1}=v,\quad \forall v\in \Gcal.\]
A quasigroup  is {\em flexible} if $u(vu)=(uv)u$ for all $u,v\in \Gcal$, and {\em Moufang} if $u(v(uw))=((uv)u)w$ for all $u,v,w\in M$.

In \cite{Kl09} we linearised these definitions to Hopf quasigroups in the same way that a Hopf algebra linearises the notion of a group.

\begin{definition}
\cite{Kl09}  A \textit{Hopf quasigroup}  is a possibly nonassociative but unital algebra $H$ equipped with algebra homomorphisms $\Delta:H\to H\tens H$, $\eps:H\to k$ forming a coassociative coalgebra and a map $S:H\to H$ such that
\begin{equation*} m(\id\tens m)(S\tens\id\tens\id)(\Delta\tens\id)=\eps\tens \id = m(\id\tens m)(\id\tens S\tens\id)(\Delta\tens\id), \end{equation*}
\begin{equation*} m(m\tens\id)(\id\tens S\tens\id)(\id\tens\Delta)=\id\tens\eps = m(m\tens\id)(\id\tens\id\tens S)(\id\tens\Delta). \end{equation*}
\end{definition}

One can write these more explicitly as
\[ \sum Sh\o (h\t g)=\sum h\o((Sh\t) g)=\sum (g Sh\o)h\t=\sum (g h\o)Sh\t=\eps(h)g,\]
for all $h,g\in H$, where we write $\Delta h=\sum h\o\tens h\t$ as in \cite{Sw69}. For brevity, we shall omit the summation signs for the remainder of the paper. In this notation the Hopf quasigroup $H$ is called \textit{flexible} if 
\[ h\o(gh\t)=(h\o g)h\t \quad \forall h,g\in H, \]
and \textit{Moufang} if 
\[\sum h\o(g(h\t f))=\sum ((h\o g)h\t)f \quad \forall h,g,f\in H.\]

In \cite{Kl09} it was proven that, as for Hopf algebras, the antipode $S$ is antimultiplicative and anticomultiplicative, that is for all $h,g\in H$
\[ S(hg)=(Sg)S(h),\quad\quad \Delta(Sh)=Sh\t\tens Sh\o.\]

\begin{lemma}
\cite{Kl09} If $H$ is a cocommutative flexible Hopf quasigroup, then $S^2=\id$ and for all $h,g\in H$,
\[ h\o(gSh\t)=(h\o g)Sh\t. \]
\label{H-adjoint}
\end{lemma}

A Hopf-like theory is established in \cite{Kl09}, for example, $S$ is antimulitpicative and anticomultiplicative. As in Hopf algebra theory, we can dualise this theory by reversing the arrows on each map to obtain a {\em Hopf coquasigroup}.

\begin{definition}
\cite{Kl09} A \textit{Hopf coquasigroup} is a unital associative algebra $A$ equipped with counital algebra homomorphisms $\Delta:A\to A\tens A$, $\eps:k\to A$, and linear map $S:A\to A$ such that
\[ (Sa\o)a\t\o\tens a\t\t = 1\tens a= a\o Sa\t\o\tens a\t\t, \]
\[ a\o\o\tens (Sa\o\t)a\t = a\tens 1 = a\o\o\tens a\o\t Sa\t, \]
for all $a\in A$. A Hopf coquasigroup is \textit{flexible} if
\[ a\o a\t\t \tens a\t\o = a\o\o a\t \tens a\o\t \quad \forall a\in A,\]
and \textit{Moufang} if
\[ a\o a\t\t\o \tens a\t\o \tens a\t\t\t = a\o\o\o a\o\t \tens a\o\o\t \tens a\t \quad \forall a\in A.\]
\end{definition}

\begin{lemma}
\cite{Kl09} If $A$ is a commutative flexible Hopf coquasigroup, then $S^2=\id$ and for all $a\in A$,
\[ a\o Sa\t\t \tens a\t\o = a\o\o Sa\t\tens a\o\t. \]
\label{A-flexible}
\end{lemma}

Actions and coactions of Hopf (co)quasigroups are also considered in \cite{Kl09}. In \cite{Br09} the author extends the theory to defining Hopf modules of Hopf (co)quasigroups and proves that the category of left and right $H$-Hopf modules are equivalent to the category of vector spaces.

\begin{definition}
\cite{Br09} Let $H$ be a Hopf quasigroup. A vector space $M$ is called  a \textit{right $H$-Hopf module} if there are maps $\alpha:M\tens H\to M$, $m\tens h\mapsto m\ra h$, and $\rho:M\to M\tens H$, $m\mapsto m\bn\tens m\bo$, such that for all $m\in M,h\in H$
\begin{equation} (m\ra h\o)\ra Sh\t=\eps(h)m=(m\ra Sh\o)\ra h\t, \quad\quad m\ra 1 = m, \label{H-quasiaction} \end{equation}
\begin{equation*} m\bn\bn\tens m\bn\bo\tens m\bo=m\bn\tens m\bo\o\tens m\bo\t, \quad\quad m\bn\eps(m\bo) = m, \end{equation*}
\begin{equation} (m\ra h)\bn\tens (m\ra h)\bo=m\bn\ra h\o\tens m\bo h\t. \label{H-compatability} \end{equation}
\end{definition}

Left $H$-Hopf modules are similarly defined. Let $H$ be a Hopf quasigroup and $M$ be a right $H$-Hopf module. Denote the set of $H$-coinvariants by $M^{coH}$, that is
\[ M^{coH}=\{ m\in M| \rho(m)=m\tens 1\}. \]
Then $M^{coH}\tens H$ is a right $H$-Hopf module with action $(m\tens h)\ra g = m\tens hg$ and coaction $m\tens h\mapsto m\tens h\o\tens h\t$.

\begin{theorem}
\cite{Br09}
If $H$ be a Hopf quasigroup and $M$ is a right $H$-Hopf module, then $M\iso M^{coH}\tens H$ as right $H$-Hopf modules.
\label{quasiHopf-module}
\end{theorem}
\proof
\cite{Br09} The proof involves constructing an isomorphism $\sig:M^{coH}\tens H\to M$ defined by $m\tens h \mapsto m\ra h$, with inverse $\sig^{-1}:M\to M^{coH}\tens H$ sending $m\mapsto m\bn\bn\ra Sm\bn\bo\tens m\bo$.
\endproof

The case of the universal differential calculus $M=\Omega^1_{univ} \subset H\tens H$ when $H$ is a Hopf coquasigroup is given in \cite{Kl09} and the proof of the required isomorphisms in the general case are similar,

\section{Integrals on Hopf quasigroups}

Let $H$ be a finite dimensional Hopf quasigroup and $H^*$ be the dual space with natural Hopf coquasigroup structure given by
\[ \lb \phi\psi,h\rb = \lb \phi,h\o \rb \lb \psi,h\t\rb, \]
\[ \lb \Delta \phi,h\tens g \rb = \lb \phi,hg\rb, \]
\[ \lb 1,h\rb=\eps(h),\]
\[ \eps (\phi) = \lb \phi,1\rb, \]
\[ \lb S\phi,h\rb=\lb \phi,Sh\rb,\]
for all $h,g\in H$, $\phi,\psi\in H^*$.

\begin{definition}
Let $H$ be a Hopf quasigroup. A \textit{left integral on $H$} is an element $\smallint \in H^*$ such that
\begin{equation} h\o\smallint (h\t) = \smallint (h).1 \quad \forall \,h\in H.\label{H-integral}\end{equation}
A left integral $\smallint$ is {\em normalised} if $\smallint(1)=1$. We denote the space of left integrals on $H$ by $\Ical^H_l$.
\end{definition}

Right integrals on $H$ are defined similarly. We find that an integral on $H$ satisfies many of the same identities as integrals on Hopf algebras.

\begin{lemma}
If $H$ is a Hopf quasigroup, then $\smallint\in H^*$ is a left integral on $H$ iff for all $\varphi\in H^*$
\[ \varphi\smallint= \eps(\varphi)\smallint.\]
\label{H-integralin}
\end{lemma}
\proof
Let $h\in H,\phi\in H^*$ then,
\[ \lb \varphi\smallint,h\rb=\lb\varphi,h\o\rb\lb\smallint,h\t\rb=\lb\varphi,h\o\smallint (h\t)\rb\]
and,
\[ \lb\eps(\varphi)\smallint,h\rb=\lb\varphi,1\rb\lb\smallint,h\rb=\lb\varphi,\smallint (h)\rb\]
Hence, $h\o\smallint (h\t) = \smallint (h).1 \iff \varphi\smallint = \eps(\varphi)\smallint$.
\endproof

Following the terminology in classical Hopf algebra theory, elements $x\in H$ such that $hx=\eps(h)x$ for all $h\in H$ are called {\em left integrals in $H$}. So by the above lemma we have that $\smallint$ is a left integral {\em on} $H$ iff $\smallint$ is a left integral {\em in} $H^*$.

\begin{lemma}
If $\smallint$ is a left integral on a Hopf quasigroup $H$ then,
\begin{enumerate}
\item $S\smallint$ is a right integral on $H$,\\
\item if $\eps(\smallint)=1$, then $\smallint =S\smallint$.
\end{enumerate}
\label{H-left-right}
\end{lemma}
\proof
For all $h\in H$ we have,
\begin{eqnarray*}
\lb S\smallint,h\o\rb h\t	&	=	&	\lb \smallint,Sh\o\rb h\t\\
	&	=	& \lb\smallint, (Sh\o)\t\rb (Sh\o)\o h\t \quad\text{by (\ref{H-integral})}\\
	&	=	& \lb\smallint, Sh\o\o\rb (Sh\o\t)h\t\\
	&	=	& \lb\smallint, Sh\rb.1\quad\text{by coassociativity}\\
	&	=	& \lb S\smallint, h\rb.1
\end{eqnarray*}
So $S\smallint$ is a right integral, which proves $(1)$. For (2), under our assumption, $\eps(\smallint)=1$ hence also $\eps(S\smallint)=1$. Therefore,
\[ \smallint=\eps(S\smallint)\smallint=(S\smallint)\smallint \quad\text{by Lemma \ref{H-integralin}}. \]
Since $S\smallint$ is a right integral, it satisfies $(S\smallint)\varphi=\eps(\varphi)S\smallint$ for all $\varphi\in H^*$ by a right integral version of Lemma \ref{H-integralin}, and in particular, $(S\smallint)\smallint = \eps(\smallint) S\smallint$. So, by our assumption, this implies that,
\[ S\smallint = \eps(\smallint) S\smallint = (S\smallint)\smallint\]
Hence, $\smallint=S\smallint$, as required.
\endproof

\begin{lemma}
If $H$ is a Hopf quasigroup and $\smallint$ is a left integral on $H$, then  for all $h,g\in H$,
\begin{enumerate}
\item $h\o\smallint(h\t Sg) = g\t\smallint(hSg\o)$,\\
\item $h\o\smallint(gh\t)=Sg\o\smallint(g\t h)$.
\end{enumerate}
\label{H-identities}
\end{lemma}
\proof
To prove $(1)$, let $h,g\in H$; by the axioms of a Hopf quasigroup,
\[ h\o\smallint(h\t Sg)=h\o\eps(g\t)\smallint(h\t Sg\o)=(h\o Sg\t\o)g\t\t\smallint(h\t Sg\o),\]
then by coassociativity and antimultiplicity of the antipode, this is equal to
\[ (h\o (Sg\o)\o)g\t\smallint(h\t(Sg\o)\t)=(h Sg\o)\o g\t\smallint(hSg\o)\t=g\t\smallint(hSg\o),\]
where the final equality is from (\ref{H-integral}).
The second identity is similar.
\endproof

\begin{theorem}
If $H$ is a finite-dimensional Hopf quasigroup then $dim(\Ical_l^H)=1$, i.e. a left integral exists and is unique up to scale.
\label{H-existence}
\end{theorem}

\proof
It was shown in \cite{Kl09} that since $H$ is finite-dimensional, $H^*$ is a Hopf coquasigroup. On $H^*$ we have a right $H$-action and right $H$-coaction given by
\begin{equation} \varphi\rharp h = \varphi\o\lb \varphi\t,Sh\rb, \label{H-action} \end{equation}
\begin{equation} \varphi\bn\tens \varphi\bo=\sum_i f^i\varphi\tens e_i, \label{H-coaction} \end{equation}
for all $h\in H$ and $\varphi\in H^*$, where $\{e_i\}$ is a basis of $H$ and $\{f^i\}$ a dual basis.
We will prove that $H^*$ is a right $H$-Hopf module with this action and coaction. Let $h,x\in H$ and $\varphi\in H^*$ then,
\begin{eqnarray*}
\lb (\varphi\rharp h\o)\rharp Sh\t,x\rb	&	=	&	\lb\varphi\rharp h\o,xS^2h\t\rb\\
	&	=	&	\lb\varphi,(xS^2h\t)Sh\o\rb\\
	&	=	&	\lb\varphi,S(h\o((Sh\t)S^{-1}x))\rb\\
	&	=	&	\lb\eps(h)\varphi,x\rb\quad\text{by the Hopf quasigroup axiom.}
\end{eqnarray*}
Hence $(\varphi\rharp h\o)\rharp Sh\t$, and similarly, $(\varphi\rharp Sh\o)\rharp h\t=\eps(h)\varphi$ for all $h\in H,\varphi\in H^*$. In order to check that (\ref{H-coaction}) is a coassociative coaction, we need $f^i\tens e_i\o\tens e_j\t=f^if^j\tens e_i\tens e_j$, which is easily verified by evaluating against general elements.
\begin{eqnarray*}
\varphi\bn\tens\varphi\bo\o\tens\varphi\bo\t	&	=	&	f^i\varphi\tens e_i\o\tens e_i\t\\
	&	=	&	(f^if^j)\varphi\tens e_i\tens e_j\\
	&	=	&	f^i(f^j\varphi)\tens e_i\tens e_j\quad\text{since $H^*$ is associative}\\
	&	=	&	f^i\varphi\bn\tens e_i\tens \varphi\bo\\
	&	=	& \varphi\bn\bn\tens \varphi\bn\bo\tens \varphi\bo.
\end{eqnarray*}
It is clear that $\varphi\bn\eps(\varphi\bo)=\varphi$, hence (\ref{H-coaction}) defines an $H$-coaction on $H^*$. Finally, we check the compatibility condition, for which we require another identity on the dual bases, $f^i\o\tens f^i\t\tens e_i=f^i\tens f^j\tens e_ie_j$.
\begin{eqnarray*}
\varphi\bn\rharp h\o\tens \varphi\bo h\t	&	=	&	(f^i\varphi)\rharp h\o\tens e_ih\t\\
	&	=	&	f^i\o\varphi\o\tens e_ih\t\lb f^i\t\varphi\t,Sh\o\rb\\
	&	=	& f^i\o\varphi\o\tens e_ih\t \lb f^i\t,Sh\o\t\rb\lb\varphi\t,Sh\o\o\rb\\
	&	=	&	f^i\varphi\o\tens (e_ie_j)h\t\t\lb f^j,Sh\t\o\rb\lb\varphi\t,Sh\o\rb\\
	&	=	&	f^i\varphi\o\tens (e_iSh\t\o)h\t\t\lb\varphi\t,Sh\o\rb\\
	&	=	&	f^i\varphi\o\tens e_i\lb\varphi\t,Sh\rb\quad\text{by the Hopf quasigroup axiom}\\
	&	=	& f^i(\varphi\rharp h)\tens e_i\\
	&	=	&	(\varphi\rharp h)\bn\tens (\varphi\rharp h)\bo
\end{eqnarray*}

By Theorem \ref{quasiHopf-module}, we deduce that $H^*\iso H^{*coH}\tens H$ as right $H$-Hopf modules, but from these definitions clearly
\[ H^{*coH} = \{ \varphi\in H^*| h\o\lb\varphi,h\t\rb = \lb \varphi,h\rb.1 \,\forall\, h\in H\} = \Ical_l^H. \]
Since, $\dim H=\dim H^*$, we can conclude that $\dim \Ical_l^H=1$.

\begin{corollary}
If H is a finite dimensional Hopf quasigroup then the antipode S is bijective.
\label{H-bijective}
\end{corollary}
\proof
This follows the proof for Hopf algebras \cite{Sw69}. Let $h\in \ker S$ and $\smallint$ be a non-zero left integral on $H$. Then under the isomorphism $\sig$ in Theorem \ref{quasiHopf-module} we have
\[ \sig(\smallint\tens h) = \smallint\rharp h = \smallint\o\lb\smallint\t,Sh\rb = \smallint\o\lb\smallint\t,0\rb = 0.\]
$\sig$ is an isomorphism hence $\smallint\tens h=0$, but $\smallint$ is non-zero, so $h=0$. Thus $\ker S=\emptyset$ and $S$ is injective. Since $H$ is finite dimensional, it is also bijective.
\endproof

There are analogous results for right integrals. We state them here without proof as we will refer to them later.

\begin{definition}
Let $H$ be a Hopf quasigroup. $\smallint_R\in H^*$ is a \textit{right integral on $H$} if for all $h\in H$,
\[ (\smallint{}_R h\o)h\t=\smallint{}_R(h).1 .\]
The space of right integrals on $H$ is denoted $\Ical_r^H$.
\end{definition}

\begin{proposition}
If $H$ is a Hopf quasigroup and $\smallint_R$ is a right integral on $H$ then,
\begin{enumerate}
\item $\smallint_R$ is a right integral on $H$ iff for all $\varphi\in H^*$, $\smallint_R\varphi=\eps(\varphi)\smallint_R$,\\
\item $\smallint_R((Sg)h\o)h\t=\smallint_R((Sg\t)h)g\o$ for all $h,g\in H$,\\
\item $\smallint_R(g\o h)g\t=\smallint_R(gh\o)Sh\t$ for all $h,g\in H$.
\end{enumerate}
\label{H-rightidentities}
\end{proposition}

\begin{theorem}
If $H$ is a finite-dimensional Hopf quasigroup then $dim(\Ical_r^H)=1$, i.e. a right integral exists and is unique up to scale.
\end{theorem}

\section{Integrals on Hopf coquasigroups}

We can similarly define integrals on a Hopf coquasigroup $A$, which satisfy equivalent properties as those for integrals on a Hopf quasigroup.

\begin{definition}
A left integral on a Hopf coquasigroup is an element $\smallint\in A^*$ such that $a\o\smallint (a\t)=\smallint(a).1$ for all $a\in A$.
\end{definition}

\begin{lemma}
Is $A$ is a Hopf coquasigroup and $\smallint\in A^*$ is a left integral on $A$ then,
\begin{enumerate}
\item $\smallint\in A^*$ is a left integral on $A$ iff for all $\varphi\in A^*$, $\varphi\smallint = \eps(\varphi)\smallint$,\\
\item $S\smallint$ is a right integral on $A$,\\
\item if $\smallint (1)=1$, then $\smallint =S\smallint$,\\
\item $a\o\smallint(a\t Sb)=b\t\smallint(aSb\o)$ for all $a,b\in A$,\\
\item $a\o\smallint(ba\t)=Sb\o\smallint(b\t a)$ for all $a,b\in A.$
\end{enumerate} 
\label{A-properties}
\end{lemma}
\proof
The proofs are analogous to those in the previous section, but we shall show the proof of (4) as an example:
\begin{eqnarray*}
a\o\smallint(a\t Sb)	&	=	&	a\o(Sb\o\t)b\t \smallint(a\t Sb\o\o)\quad\text{by the Hopf coquasigroup axiom}\\
	&	=	&	a\o (Sb\o)\o b\t\smallint(a\t (Sb\o)\t)\\
	&	=	&	(aSb\o)\o b\t\smallint(a Sb\o)\t\\
	&	=	&	b\t\smallint(a Sb\o)
\end{eqnarray*}
\endproof

In order to prove the existence of these left integrals, we first need to prove a form of the Hopf module lemma on Hopf coquasigroups.

\begin{definition}
\cite{Br09} Let $A$ be a Hopf coquasigroup. A vector space $M$ is a right $A$-Hopf module if there are maps $\alpha:M\tens A\to M,\,m\tens a\mapsto m\ra a$ and $\rho:M\to M\tens A,\, m\mapsto m\bn\tens m\bo$, such that $\alpha$ is an associative unital right-action on $M$ and
\begin{equation} m\bn\bn\tens (Sm\bn\bo)m\bo=m\tens 1=m\bn\bn\tens m\bn\bo Sm\bo, m\bn\eps(m\bo)=m \label{A-quasicoaction}\end{equation}
\begin{equation} (m\ra a)\bn\tens (m\ra a)\bo=m\bn\ra a\o\tens m\bo a\t, \label{A-compatibility}\end{equation}
for all $a\in A,m\in M$.
\end{definition}

In \cite{Br09} it was shown that if $M$ is a right $A$-Hopf module with coaction $\rho(m)=m\bn\tens m\bo$, then $M$ is a right $A$-Hopf module with the same action and the induced coaction
\[ \widehat{\rho}(m)=m\bn\bn\ra((Sm\bn\bo)m\bo\o)\tens m\bo\t. \]

\begin{theorem}
If $A$ is a Hopf coquasigroup and $M$ is a right $A$-Hopf module, then $M\iso M^{\widehat{coA}}\tens A$ as right $A$-Hopf modules, where $M^{\widehat{coA}}\tens A$ is a right $A$-Hopf module by
\[ (m\tens a)\ra b=m\tens ab, \quad\quad (m\tens a)\bn\tens (m\tens a)\bo=m\tens a\o\tens a\t,\]
for all $a,b\in A,m\in M$, and $M^{\widehat{coA}}=\{m\in M|\widehat{\rho}(m)=m\tens 1\}$. The isomorphism is given by
\[ \sig:M^{\widehat{coA}}\tens A\to M \quad\text{by} \quad m\tens a\mapsto m\ra a,\]
\[\sig^{-1}:M\to M^{\widehat{coA}}\tens A \quad\text{by}\quad m\mapsto m\bn\bn\ra Sm\bn\bo\tens m\bo.\]
\label{coquasiHopf-module}
\end{theorem}
\proof
It is straightforward to check that $\sigma$ as stated gives the required isomorphism, however here we will make use of Theorem 3.13 in \cite{Br09}, which gives an isomorphism $M\iso M^A\tens A$, where $M^A$ is the set of $A$-invariants defined by the coequalizer

\centerline{
\xymatrix{
{M\tens A} \ar@<.7ex>[r]^{\alpha} \ar@<-.7ex>[r]_{\id\tens\eps}
	&	M  \ar[r]^{\pi_M}
		&	{M^A}
}
}
so that $\pi_M\circ\alpha=\pi_M(\id\tens\eps)$. As noted in \cite{Br09}, for Hopf algebras, there exists an isomorphism $M^A\iso M^{coA}$. We now show that the same maps also provide an isomorphism $M^A\iso M^{\widehat{coA}}$ when $A$ is a Hopf coquasigroup; define $\omega:M^{\widehat{coA}}\to M^A$ by $m\mapsto \pi_M(m)$ for $m\in M^{\widehat{coA}}$, and $\omega^{-1}:M^A\to M^{\widehat{coA}}$ by $\pi_M(m)\mapsto m\bn\ra Sm\bo$ for all $m\in M$.

First we have to check that $\omega^{-1}$ is well-defined,
\begin{eqnarray*}
\lefteqn{\widehat{\rho}(\omega^{-1}(\pi_M(m)))=}\\
	&	=	&	\widehat{\rho}(m\bn\ra Sm\bo)\\
	&	=	&	(m\bn\ra Sm\bo)\bn\bn\ra((S(m\bn\ra Sm\bo)\bn\bo)(m\bn\ra Sm\bo)\bo\o)\tens (m\bn\ra Sm\bo)\bo\t\\
	&	=	&	(m\bn\bn\bn\ra(Sm\bo)\o\o)\ra(S(m\bn\bn\bo(Sm)\o\t)(m\bn\bo(Sm\bo)\t)\o)\\
	&		&	\quad\tens (m\bn\bo (Sm\bo)\t)\t \quad\text{by the compatibility condition (\ref{A-compatibility}),}\\
	&	=	&	m\bn\bn\bn\ra((Sm\bo)\o\o S((Sm\bo)\o\t)(Sm\bn\bn\bo)(m\bn\bo(Sm\bo)\t)\o)\\
	&		&	\quad\tens (m\bn\bo(Sm\bo)\t)\t\\
	&	=	&	m\bn\bn\bn\ra((Sm\bn\bn\bo)(m\bn\bo Sm\bo)\o)\tens(m\bn\bo Sm\bo)\t\\
	&		&	\quad\quad\text{by the antipode property on $(Sm\bo)\o$}\\
	&	=	&	m\bn\ra Sm\bo\tens 1 \quad\text{by (\ref{A-quasicoaction}),}\\
	&	=	& \omega^{-1}(\pi_M(m))\tens 1.
\end{eqnarray*}
Hence, $\omega^{-1}(\pi_M(m))$ lies in $M^{\widehat{coA}}$ so is well-defined. Finally, we check that these maps are mutually inverse; for all $m\in M$,
\begin{eqnarray*}
\omega(\omega^{-1}(\pi_M(m)))	&	=	&	\omega(m\bn\ra Sm\bo)\\
	&	=	&	\pi_M(m\bn\ra Sm\bo)\\
	&	=	&	\pi_M(m\bn\eps(Sm\bo))\\
	&	=	&	\pi_M(m)\quad\text{by (\ref{A-quasicoaction}),}
\end{eqnarray*}
and for all $m\in M^{\widehat{coA}}$,
\begin{eqnarray*}
\omega^{-1}(\omega(m))	&	=	&	\omega^{-1}(\pi_M(m))\\
	&	=	&	m\bn\ra Sm\bo\\
	&	=	&	m\bn\bn\ra Sm\bn\bo \eps(m\bo) \quad\text{by (\ref{A-quasicoaction}),}\\
	&	=	&	m\bn\bn\ra (Sm\bn\bo)m\bo\o Sm\bo\t\\
	&	=	&	(m\bn\bn\ra (Sm\bn\bo)m\bo\o)\ra Sm\bo\t\\
	&	=	&	m\ra S1\quad\text{since $m\in M^{\widehat{coA}}$,}\\
	&	=	&	m.
\end{eqnarray*}
\endproof

\begin{theorem}
If $A$ is a finite dimensional, commutative, flexible Hopf coquasigroup, then $dim(\Ical_l^A)=1$, i.e. a left integral exists and is unique up to scale.
\label{A-existence}
\end{theorem}
\proof
Since $A$ is finite dimensional, $A^*$ is a Hopf quasigroup and hence coassociative. On $A^*$ we have a right $A$-action and a right $A$-coaction given as in Theorem \ref{H-existence} by
\begin{equation} \varphi\rharp a=\varphi\o\lb\varphi\t,Sa\rb,\label{A-action}\end{equation}
\begin{equation} \varphi\bn\tens\varphi\bo=f^i\varphi\tens e_i.\label{A-coaction}\end{equation}
for $a\in A,\varphi\in A^*$, where $\{e_i\}$ is a basis for $A$ and $\{f^i\}$ is a dual basis. The proof is similar to the case when $A$ is a Hopf quasigroup; let $\varphi\in A^*$ and $a,b\in A$ then
\begin{eqnarray*}
(\varphi\rharp a)\rharp b	&	=	&	\varphi\o\rharp b\lb\varphi\t,Sa\rb\\
	&	=	&	\varphi\o\o\lb\varphi,Sb\rb\lb\varphi\t,Sa\rb\\
	&	=	&	\varphi\o\lb\varphi\t\o,Sb\rb\lb\varphi\t\t,Sa\rb\quad\text{by coassociativity in $A^*$,}\\
	&	=	&	\varphi\o\lb\varphi\t,S(ab)\rb\\
	&	=	&	\varphi\rharp(ab)
\end{eqnarray*}
and clearly, $\varphi\rharp 1=\varphi$, so (\ref{A-action} is a associative unital right $A$-action. For the coaction we evaluate against $a\in A$ in the first factor and find,
\begin{eqnarray*}
\lb\varphi\bn\bn,a\rb(S\varphi\bn\bo)\varphi\bo	&	=	&	\lb (f^i\varphi)\bn,a\rb(S(f^i\varphi)\bo)e_i\\
	&	=	&	\lb f^j(f^i\varphi),a\rb(Se_j)e_i\\
	&	=	&	\lb f^j,a\o\rb\lb f^i,a\t\o\rb\lb\varphi,a\t\t\rb(Se_j)e_i\\
	&	=	&	\lb\varphi,a\t\t\rb(Sa\o)a\t\o\\
	&	=	&	\lb\varphi,a\rb.1
\end{eqnarray*}
This holds for all $a\in A$, hence, $\varphi\bn\bn\tens (S\varphi\bn\bo)\varphi\bo=\varphi\tens 1$. Similarly, $\varphi\bn\bn\tens \varphi\bn\bo S\varphi\bo=\varphi\tens 1$. Finally, we have to check the compatibility condition holds; let $a,b\in A$ and $\varphi\in A^*$ then,
\begin{eqnarray*}
\lb\varphi\bn\rharp a\o,b\rb\varphi\bo a\t	&	=	&	\lb (f^i\varphi)\rharp a\o,b\rb e_ia\t\\
	&	=	&	\lb f^i\o\varphi\o,b\rb \lb f^i\t\varphi\t,Sa\o\rb e_ia\t\\
	&	=	&	\lb f^i\varphi\o,b\rb\lb f^j,Sa\o\t\rb\lb\varphi\t,Sa\o\o\rb e_ie_ja\t\\
	&	=	&	\lb f^i\varphi\o,b\rb\lb\varphi\t,Sa\o\o\rb e_i(Sa\o\t)a\t\\
	&	=	&	\lb f^i\varphi\o,b\rb\lb\varphi\t,Sa\rb e_i\quad\text{by the Hopf coquasigroup axiom,}\\
	&	=	&	\lb f^i(\varphi\rharp a),b\rb e_i\\
	&	=	&	\lb (\varphi\rharp a)\bn,b\rb(\varphi\rharp a)\bo.
\end{eqnarray*}
Hence, $(\varphi\rharp a)\bn\tens(\varphi\rharp a)\bo=\varphi\bn\rharp a\o\tens \varphi\bo a\t$ for all $a\in A,\varphi\in A^*$, and $A^*$ is a right $A$-Hopf module. We note that for this right coaction $\rho$, $\widehat{\rho}=\rho$; let $a\in A$ and $\varphi\in A^*$ then,
\begin{eqnarray*}
\lb\varphi\bn\bn\rharp((S\varphi\bn\bo)\varphi\bo\o),a\rb\varphi\bo\t	&	=	&	\lb\varphi\bn\bn,aS((S\varphi\bn\bo)\varphi\bo\o)\rb\varphi\bo\t\\
	&	=	&	\lb\varphi,a\t\t S((Sa\o)a\t\o\o)\rb a\t\o\t\\
	&	=	&	\lb\varphi,S((Sa\o)a\t\o\o S^{-1}a\t\t)\rb a\t\o\t\\
	&	=	&	\lb\varphi,S((Sa\o)a\t\o\o Sa\t\t)\rb a\t\o\t\\
	&		&	\quad\quad\text{since $S^2=\id$ Lemma \ref{A-flexible}}\\
	&	=	&	\lb\varphi,S((Sa\o)a\t\o Sa\t\t\t)\rb a\t\t\o\\
	&		&	\quad\quad\text{by Lemma \ref{A-flexible} on $a\t$,}\\
	&	=	&	\lb\varphi,S^2a\t\rb a\o\quad\text{by the Hopf coquasigroup axiom,}\\
	&	=	&	\lb\varphi,a\t\rb a\o\\
	&	=	&	\lb\varphi\bn,a\rb \varphi\bo.
\end{eqnarray*}
By Theorem \ref{coquasiHopf-module}, we deduce that $A^*\iso A^{*\widehat{coA}}\tens A=A^{*coA}\tens A$, and as in Lemma \ref{H-existence},
\[ A^{*coA}=\{\varphi\in A^*|a\o\lb\varphi,a\t\rb=\lb\varphi,a\rb.1 \,\forall \, a\in A\}=\Ical_l^A.\]
\endproof

\section{Fourier transformations on Hopf quasigroups}

In this section $H$ will be a finite dimensional Hopf quasigroup, and $H^*$ its dual space with the structure of a Hopf coquasigroup. We will follow the treatment in \cite{Ma95} of Fourier transformations on Hopf algebras.

\begin{definition}
Let $\smallint$ be a left integral on $H$. A {\em Fourier transformation} on $H$ is a linear map $F:H\to H^*$ defined by
\begin{equation} F(h) = \smallint\rharp h, \label{H-FT}\end{equation}
for all $h\in H$, where the right action of $H$ on $H^*$ is given by (\ref{H-action}).
\end{definition}

\begin{lemma}
If $\smallint$ is a left integral on $H$ and $F$ is a Fourier transform on $H$, then for all $h\in H$, $\varphi\in H^*$,
\begin{enumerate}
	\item $F(h\o Sh\t) = F(h\o)\rharp Sh\t$,\\
	\item $F((Sh\o) h\t) = F(Sh\o)\rharp h\t$,\\
	\item $F(\varphi\lharp h) = \varphi F(h)$,\\
	\item $\rho(F(h))=(F\tens\id)\Delta(h)$,\\
	\item $H$ cocommutative and flexible $\implies$ $\lb F(h\o g),h\t\rb=\lb F(h\o)\rharp g,h\t\rb,$
\end{enumerate}
where the left coaction is given by (\ref{H-coaction}) and the left action of $H^*$ on $H$ is given by $\varphi\lharp h=h\o\lb\varphi,h\t\rb$ for all $h\in H$ and $\varphi\in H^*$.
\label{H-FT-properties}
\end{lemma}
\proof
To prove $(1)$ we evaluate against a general element $x\in H$,
\begin{eqnarray*}
\lb F(h\o Sh\t),x\rb	&	=	&	\lb \smallint, xS(h\o Sh\t)\rb\\
	&	=	& \lb \smallint, x(S(Sh\t)Sh\o)\rb\\
	&	=	& \lb \smallint, x(S((Sh)\o)(Sh)\t)\rb\\
	&	=	& \lb \smallint, (xS(Sh)\o)(Sh)\t\rb\quad\text{by the Hopf quasigroup axiom,}\\
	&	=	&	\lb \smallint, (xS(Sh\t))Sh\o\rb\\
	&	=	& \lb F(h\o),xS(Sh\t)\rb\\
	&	=	& \lb F(h\o)\rharp Sh\t,x\rb\quad\text{by (\ref{H-action}).}
\end{eqnarray*}
Therefore $F(h\o Sh\t) = F(h\o)\rharp Sh\t$. Similarly we can prove $(2)$. To prove $(3)$, let $x\in H$ then,
\begin{eqnarray*}
\lb F(\varphi\lharp h),x\rb	&	=	&	\lb F(h\o),x\rb\lb\varphi,h\t\rb\\
	&	=	&	\lb\smallint,xSh\o\rb\lb\varphi,h\t\rb\\
	&	=	&	\lb\smallint,x\t Sh\rb\lb\varphi,x\o\rb\quad\text{by Lemma \ref{H-identities} (1),}\\
	&	=	&	\lb F(h),x\t\rb\lb\varphi,x\o\rb\\
	&	=	&	\lb\varphi F(h),x\rb.
\end{eqnarray*}
This holds for all $x,h\in H$, hence $F(\varphi\lharp h)=\varphi F(h)$ for all $h\in H,\varphi\in H^*$. To prove $(4)$,
\begin{eqnarray*}
\rho(F(h))	&	=	& f^iF(h)\tens e_i\\
	&	=	&	F(f^i\lharp h)\tens e_i\quad\text{by (3),}\\
	&	=	&	F(h\o)\tens e_i\lb f^i,h\t\rb\\
	&	=	&	F(h\o)\tens h\t\\
	&	=	&	(F\tens\id)\Delta(h)
\end{eqnarray*}
Finally if $H$ is cocommutative and flexible we find,
\begin{eqnarray*}
\lb F(h\o)\rharp g,h\t\rb	&	=	&	\lb F(h\o),h\t Sg\rb\\
	&	=	&	\lb\smallint,(h\t Sg)Sh\o\rb\\
	&	=	&	\lb\smallint,(h\o Sg)Sh\t\rb\quad\text{since $H$ is cocommutative,}\\
	&	=	&	\lb\smallint,h\o(Sg Sh\t)\rb\quad\text{by Lemma \ref{H-adjoint},}\\
	&	=	&	\lb\smallint, h\o S(h\t g)\rb\\
	&	=	&	\lb \smallint,h\t S(h\o g)\rb\\
	&	=	&	\lb F(h\o g),h\t\rb
\end{eqnarray*}
which is the required identity $(5)$.
\endproof

\begin{lemma}
If $H$ is a Hopf quasigroup with Fourier transform $F$ defined by (\ref{H-FT}), then for all $g,h\in H$,
\[ F(g)F(h) = F(F(g)\lharp h). \]
\label{H-FT-product}
\end{lemma}
\proof
Let $x\in H$, then
\begin{eqnarray*}
\lb F(g)F(h),x\rb	&	=	&	\lb F(g),x\o\rb\lb F(h),x\t\rb\\
	&	=	& \lb F(g),x\o\rb\lb\smallint,x\t Sh\rb\\
	&	=	&	\lb F(g),h\t\rb\lb\smallint,xSh\o\rb\quad\text{by Lemma \ref{H-identities} (1),}\\
	&	=	&	\lb \smallint,x S(F(g)\lharp h)\rb\\
	&	=	&	\lb F(F(g)\lharp h),x\rb.
\end{eqnarray*}
Since this holds for all $x\in H$, $F(g)F(h) = F(F(g)\lharp h)$.
\endproof

Let $\smallint$ be a left integral on $H$ and $\smallint^*_R$ be a right integral on $H^*$ and set $\mu =\lb \smallint,\smallint^*_R\rb$. Define the map $F^{-1}:H^*\to H$ by
\begin{equation} F^{-1}(\varphi)= \frac{1}{\mu}\,\,(\varphi\lharp\smallint{}^*_R),\label{H-FT-inverse}\end{equation}
for all $\varphi\in H^*$. We show that $F$ and $F^{-1}$ are mutually inverse.
\begin{eqnarray*}
F^{-1}(F(h))	&	=	&	\frac{1}{\mu} (\smallint\rharp h)\lharp\smallint{}^*_R\\
	&	=	&	\frac{1}{\mu}\smallint{}^*_R\o\lb\smallint\rharp h,\smallint{}^*_R\t\rb\\
	&	=	&	\frac{1}{\mu}\smallint{}^*_R\o\lb\smallint,\smallint{}^*_R\t Sh\rb\\
	&	=	&	\frac{1}{\mu}h\t\lb\smallint,\smallint{}^*_R Sh\o\rb\quad\text{by Lemma \ref{H-identities} (1),}\\
	&	=	&	\frac{1}{\mu}h\t\lb\smallint,\eps(Sh\o)\smallint{}^*_R\rb\quad\text{since $\smallint{}^*_R$ is a right integral,}\\
	&	=	&	\frac{1}{\mu} h\lb\smallint,\smallint{}^*_R\rb\\
	&	=	&	h.
\end{eqnarray*}
Also,
\begin{eqnarray*}
F(F^{-1}(\varphi))	&	=	& \frac{1}{\mu}\smallint\rharp(\varphi\lharp\smallint{}^*_R)\\
	&	=	&	\frac{1}{\mu}\smallint\o\lb\smallint\t,S(\varphi\lharp\smallint{}^*_R)\rb\\
	&	=	&	\frac{1}{\mu}\smallint\o\lb S\smallint\t,\varphi\lharp\smallint{}^*_R\rb\\
	&	=	&	\frac{1}{\mu}\smallint\o\lb(S\smallint\t)\varphi,\smallint{}^*_R\rb\\
	&	=	&	\frac{1}{\mu}\varphi\t\lb(S\smallint)\varphi\o,\smallint{}^*_R\rb\quad\text{by Lemma \ref{H-identities} (2),}\\
	&	=	&	\frac{1}{\mu}\varphi\t\lb \eps(\varphi\o)S\smallint,\smallint{}^*_R\rb\\
	&		&	\quad\quad\text{since $S\smallint$ is a right integral by Lemma \ref{H-left-right} (1),}\\
	&	=	&	\frac{1}{\mu}\varphi\lb S\smallint,\smallint{}^*_R\rb
\end{eqnarray*}
Now, $F^{-1}\circ F=\id$, so $F^{-1}(F(F^{-1}(\varphi)))=F^{-1}(\varphi)$, but also, by the above calculation,
\[ F^{-1}(F(F^{-1}(\varphi)))=F^{-1}(\frac{1}{\mu}\lb S\smallint,\smallint{}^*_R\rb\varphi)=\frac{1}{\mu}\lb S\smallint,\smallint{}^*_R\rb F^{-1}(\varphi).\]
So we have
\[ F^{-1}(\varphi)=\frac{1}{\mu}\lb S\smallint,\smallint{}^*_R\rb F^{-1}(\varphi),\]
which implies $\lb S\smallint,\smallint^*_R\rb=\mu=\lb\smallint,\smallint^*_R\rb$, and hence $F(F^{-1}(\varphi))=\varphi$ and the two maps are mutually inverse.

\begin{definition}
For $g,h\in H$ and $\smallint\in H^*$ a left integral on $H$, define the {\em convolution product} $g*h\in H$ by
\[ g*h = h\o\lb\smallint,h\t Sg\rb. \]
\end{definition}

\begin{proposition}
If $H$ is a finite dimensional Hopf-quasigroup with non-zero left integral $\smallint$, then the Fourier transform $F:H\to H^*$ maps the convolution product in $H$ to the product in $H^*$, that is,
\[F(g*h)=F(g)F(h), \quad\quad \forall h,g\in H.\]
\end{proposition}
\proof
Let $g,h,x\in H$ then,
\begin{eqnarray*}
\lb F(g)F(h),x\rb	&	=	& \lb F(F(g)\lharp h),x\rb\quad\text{by Lemma \ref{H-FT-product},}\\
	&	=	&	\lb F(h\o),x\rb\lb F(g),h\t\rb\\
	&	=	&	\lb F(h\o),x\rb\lb\smallint,h\t Sg\rb\\
	&	=	&	\lb F(h\o\lb\smallint,h\t Sg\rb),x\rb\\
	&	=	&	\lb F(g*h),x\rb.
\end{eqnarray*}
This holds for all $x\in H$, hence, $F(g)F(h)=F(g*h)$ for all $g,h\in H$.
\endproof

\begin{lemma}
If $\smallint$ is a left integral on $H$ such that $\smallint(hg)=\smallint(gS^2h)$ for all $h,g\in H$, then
\[ S(g*h)=S(h)*S(g),\]
for all $g,h\in H$.
\label{H-FT-antipode}
\end{lemma}
\proof
Let $h,g\in H$ then,
\begin{eqnarray*}
S(g*h)	&	=	&	Sh\o\lb\smallint,h\t Sg\rb\\
	&	=	&	Sh\o\lb\smallint, (Sg)S^2h\t\rb\quad\text{by assumption,}\\
	&	=	&	Sh\o\lb S\smallint, (Sh\t)g\rb\\
	&	=	&	Sg\t\lb S\smallint, (Sh)g\o\rb\quad\text{by Proposition \ref{H-rightidentities} (2),}\\
	&	=	&	Sg\t\lb\smallint, (Sg\o)S^2h\rb\\
	&	=	&	(Sg)\o\lb\smallint, (Sg)\t S^2h\rb\\
	&	=	&	(Sh)*(Sg).
\end{eqnarray*}
\endproof

\section{Fourier transformations on Hopf coquasigroups}

Let $A$ be a finite dimensional, commutative, flexible Hopf coquasigroup. In section 4 we showed that there exists a unique left integral $\smallint$ on $A$. We define a Fourier transform on a Hopf coquasigroup in the same way as for Hopf quasigroups.

\begin{definition}
Let $\smallint$ be a left integral on $A$. A {\em Fourier transformation} on $A$ is a linear map $F:A\to A^*$ defined by
\begin{equation} F(a) = \smallint\rharp a, \label{A-FT}\end{equation}
for all $a\in A$, where the action of $A$ on $A^*$ is given by (\ref{A-action}). The inverse Fourier transform $F:A^*\to A$ is given by
\begin{equation} F^{-1}(\varphi)= \frac{1}{\mu}\,\,(\varphi\lharp\smallint{}^*_R),\label{A-FT-inverse}\end{equation}
for all $\varphi\in A^*$ where $\smallint{}^*_R$ is a right integral on $A^*$, $\mu=\lb\smallint,\smallint{}^*_R\rb$, and the left action $\lharp$ is as defined in Lemma \ref{H-FT-properties}.
\end{definition}

The proof that these maps are indeed mutually inverse is the same as for Hopf quasigroups as we did not require coassociativity of the Hopf quasigroup, and the Fourier transform has similar properties as before.

\begin{lemma}
If $\smallint$ is a left integral on $A$ and $F$ is a Fourier transform on $A$, then for all $a,b\in A$, $\varphi\in A^*$,
\begin{enumerate}
	\item $F(ab)=F(a)\rharp b$,\\
	\item $F(\varphi\lharp a)=\varphi F(a)$,\\
	\item $\rho(F(a))=(F\tens\id)\Delta(a)$,\\
	\item $F(a)F(b)=F(F(a)\lharp b)$,
\end{enumerate}
where the left coaction is given by (\ref{A-coaction}).
\label{A-FT-properties}
\end{lemma}
\proof 
The proofs of $(2)$, $(3)$ and $(4)$ are the same as those in Lemma \ref{H-FT-properties} so we shall omit them. The proof of $(1)$ just uses associativity of the action. Let $a,b\in A$ then
\[ F(ab)=\smallint\rharp(ab)=(\smallint\rharp a)\rharp b=F(a)\rharp b.\]
\endproof

The convolution product on a Hopf quasigroup $A$ is defined as before by
\[ a*b = b\o\lb\smallint, b\t Sa\rb,\]
for all $a,b\in A$.

\begin{proposition}
If $A$ is a finite dimensional Hopf coquasigroup with non-zero left integral $\smallint$, then the Fourier transform $F:A\to A^*$ maps the convolution product in $A$ to the product in $A^*$, that is,
\[F(a*b)=F(a)F(b), \quad\quad \forall a,b\in A.\]
\end{proposition}
\proof
The proof is as in Proposition 5.5 since we did not require coassociativity of the Hopf quasigroup. We include it here for completeness.
Let $a,b,x\in A$ then,
\begin{eqnarray*}
\lb F(a)F(b),x\rb	&	=	& \lb F(F(a)\lharp b),x\rb \quad\text{by Lemma \ref{A-FT-properties} (4),}\\
	&	=	&	\lb F(b\o),x\rb\lb F(a),b\t\rb\\
	&	=	&	\lb F(b\o),x\rb\lb\smallint, b\t Sa\rb\\
	&	=	&	\lb F(b\o\lb\smallint,b\t Sa\rb),x\rb\\
	&	=	&	\lb F(a*b),x\rb.\\
\end{eqnarray*}
This holds for all $x\in H$, hence, $F(a)F(b)=F(a*b)$ for all $a,b\in A$.
\endproof

\begin{lemma}
If $\smallint$ is a left integral on $A$ such that $\smallint(ab)=\smallint(bS^2a)$ for all $a,b\in A$, then
\[ S(a*b)=(Sb)*(Sa),\]
for all $a,b\in A$.
\end{lemma}
\proof
The proof is as in Lemma \ref{H-FT-antipode}.
\endproof

\begin{example}
Let $\Gcal$ be a finite (IP) quasigroup with identity element $e$; by definition $u^{-1}(uv)=v=(vu)u^{-1}$ for all $u,v\in\Gcal$. It was shown in \cite{Kl09} that the group algebra $k\Gcal$ is a Hopf quasigroup. Since $\Gcal$ is finite, the dual space $k[\Gcal]$ is a Hopf coquasigroup with structure shown here on the basis elements $\{\delta_s|s\in\Gcal\}$.
\[\delta_s\delta_t=\delta_{s,t}\,\delta_t\]
\[1=\sum_t\delta_t\]
\[\Delta(\delta_s)=\sum_t\delta_t\tens\delta_{t^{-1}s}\]
\[\eps(\delta_s)=\delta_{s,e}\]
\[ S(\delta_s)=\delta_{s^{-1}}\]
Then $\sum_{u\in\Gcal} u$ is a left integral for $k[\Gcal]$ since,
\[ \delta_s\o\lb\delta_s\t,\sum_{u\in\Gcal}u\rb=\sum_{t,u}\delta_t\delta_{t^{-1}s,u}=\sum_u\delta_{su^{-1}}=1=\lb\delta_s,u\rb.1\]
The Fourier transformation on the basis elements is given by
\[ F(\delta_s)=\sum_u (u\rharp\delta_s)=\sum_u u\o \delta_{s^{-1},u\t}=\sum_u u\delta_{s^{-1},u}=s^{-1},\]
and the convolution product by
\[ \delta_s*\delta_t=\sum_u\delta_t\o\lb\delta_t\t S\delta_s,u\rb=\sum_{u,x} \delta_x\lb\delta_{x^{-1}t}\delta_{s^{-1}},u\rb=\sum_x\delta_x\delta_{x^{-1}t,s^{-1}}=\delta_{ts}.\]

\end{example}

\section{Frobenius and separable Hopf (co)quasigroups}

\begin{definition}
Let $H$ be a Hopf quasigroup and $B:H\times H\to k$ be a bilinear form on $H$. $B$ is {\em inverse associative} if
\[ B(h\o,(Sh\t)g)=\eps(h)B(1,g)=B(Sh\o,h\t g),\]
\[ B(gh\o,Sh\t)=\eps(h)B(g,1)=B(gSh\o,h\t),\]
for all $g,h\in H$, 
and we say it is {\em non-degenerate} if
\[ B(h,g)=0\,\forall g\in H\implies h=0, \quad\text{and}\quad B(h,g)=0\,\forall h\in H\implies g=0.\]
\end{definition}

\begin{definition}
A Hopf quasigroup $H$ is {\em Frobenius} if there exists an inverse associative non-degenerate bilinear form on $H$.
\end{definition}

\begin{proposition}
Every finite dimensional Hopf quasigroup is Frobenius.
\end{proposition}
\proof
Let $H$ be a finite dimensional Hopf quasigroup, then by Theorem \ref{H-existence}, a left integral $\int\in H^*$ on $H$ exists and is unique up to normalization. By Theorem \ref{quasiHopf-module} we have an isomorphism $H^*\iso \int\tens H$. Recall the left action of $H$ on $H^*$ defined in Lemma \ref{H-FT-properties} by $h\lharp\varphi=\varphi\o\lb\varphi\t,h\rb$, so $\varphi\rharp h=Sh\lharp\varphi$ for all $h\in H,\varphi\in H^*$, where the right action is defined in (\ref{H-action}). Then we have
\[ H^*= \smallint\rharp H=SH\lharp \smallint=H\lharp\smallint,\]
the last equality coming from the bijectivity of the antipode from Corollary \ref{H-bijective}. Now, define a form $B:H\times H\to k$ on $H$ by $B(h,g)=\lb\int,hg\rb$ for all $h,g\in H$. Clearly $B$ is bilinear, and it is inverse associative by a simple application of the Hopf quasigroup axioms, for example,
\[ B(Sh\o,h\t g)=\lb\smallint, Sh\o(h\t g)\rb=\lb\smallint,\eps(h)g\rb=\eps(h)B(1,g).\]
It remains to show that $B$ is non-degenerate. Assume there exists $x\in H$ such that $B(x,H)=0$, then,
\[ 0=B(x,H)=\lb\smallint,xH\rb=\lb H\lharp\smallint,x\rb=\lb H^*,x\rb.\]
By the non-degeneracy of the dual pairing, this implies that $x=0$, and $B$ is left non-degenerate. Since $H$ is finite dimensional this is sufficient to prove that $B$ is non-degenerate.
\endproof

\begin{definition}
Let $A$ be a Hopf coquasigroup and $B:A\times A\to k$ be a bilinear form on $A$. We say $B$ is {\em associative} if 
\[ B(a,bc)=B(ab,c) \quad\quad \forall a,b,c\in A,\]
and we say it is {\em non-degenerate} if
\[ B(a,b)=0\,\forall b\in H\implies a=0, \quad\text{and}\quad B(a,b)=0\,\forall a\in H\implies b=0.\]
\end{definition}

\begin{definition}
A Hopf coquasigroup $A$ is {\em Frobenius} if there exists an associative non-degenerate bilinear form on $A$.
\end{definition}

\begin{proposition}
Every finite dimensional, commutative, flexible Hopf coquasigroup is Frobenius.
\end{proposition}
\proof
If $A$ is a finite dimensional, commutative, flexible Hopf coquasigroup, the by Theorem \ref{A-existence}, a left integral $\int\in A^*$ on $A$ exists and is unique up to scale. Define a form $B:A\times A\to k$ by $B(a,b)=\lb\int,ab\rb$ for $a,b\in A$. Clearly $B$ is associative since $A$ is;
\[ B(a,bc)=\lb\smallint,a(bc)\rb=\lb\smallint,(ab)c\rb=B(ab,c).\]
The proof that $B$ is non-degenerate is the same as in the proof that finite dimensional Hopf quasigroups are Frobenius.
\endproof

\begin{definition}
A Hopf (co)quasigroup $H$ is {\em separable} if there exists $\omega=\sum\omega\bo\tens\omega\bt\in H\tens H$, such that $\sum\omega\bo\omega\bt=1$ and $\sum h\omega\bo\tens\omega\bt=\sum\omega\bo\tens\omega\bt h$ for all $h\in H$.
\end{definition}

\begin{proposition}
If $H$ is a Hopf (co)quasigroup with left integral $\Lambda\in H$ {\em in} $H$, then $H$ is separable.
\end{proposition}
\proof
Let $\Lambda\in H$ be a left integral in $H$ with $\eps(\Lambda)=1$, that is $h\Lambda=\eps(h)\Lambda$ for all $h\in H$. Define $\omega=\Lambda\o\tens S\Lambda\t$. For all $h\in H$,
\[ \Delta(\Lambda)\tens h=\Delta(\eps(h\o)\Lambda)\tens h\t=\Delta(h\o\Lambda)\tens h\t,\]
so,
\begin{eqnarray*}
\Lambda\o\tens (S\Lambda\t)h	&	=	& h\o\o\Lambda\o\tens S(h\o\t\Lambda\t)h\t\\
	&	=	& h\o\o\Lambda\o\tens ((S\Lambda\t)Sh\o\t)h\t\\
	&	=	& h\Lambda\o\tens S\Lambda\t.\quad\text{by Hopf (co)quasigroup axioms}
\end{eqnarray*}
Also, $\omega\bo\omega\bt=\Lambda\o S\Lambda\t=\eps(\Lambda)=1$, and $\omega$ is the required separability element.
\endproof

\begin{corollary}
\begin{enumerate}
	\item Every finite dimensional, cocommutative, flexible Hopf quasigroup is separable,\\
	\item Every finite dimensional Hopf coquasigroup is separable.
\end{enumerate}
\end{corollary}
\proof
If a Hopf quasigroup $H$ is as in $(1)$, then by Theorem \ref{A-existence}, a left integral {\em on} $H^*$, and hence a left integral {\em in} $H$ exists. While if $H$ is a finite dimensional Hopf coquasigroup as in $(2)$, Theorem \ref{H-existence} proves that a left integral {\em in} $H$ exists. 
\endproof

\section{Semisimple Hopf coquasigroups}

\begin{definition}
A Hopf coquasigroup $A$ is semisimple iff every left $A$-module is completely reducible.
\end{definition}

\begin{proposition}
Let $A$ be a finite-dimensional Hopf coquasigroup and $\Lambda$ be a non-singular left integral in $A$. Then $A$ is semi-simple {\em iff} $\eps(\Lambda)\ne 0$.
\end{proposition}

\proof
The proof is as in \cite{La69}, although we now require the Hopf coquasigroup axioms. Suppose $\eps(\Lambda)=0$. Then $\Lambda^2=\eps(\Lambda)\Lambda=0$, so $k\Lambda$ is a nilpotent left ideal in $A$. Therefore, if $A$ is semisimple then $\eps(\Lambda)\ne 0$.

Conversely, suppose $\eps(\Lambda)\ne 0$. By uniqueness of the integral, we can assume $\eps(\Lambda)=1$. To prove $A$ is semisimple, it is sufficient to prove that every left $A$-module is completely reducible. Let $M$ be a left $A$-module (by an associative, unital action labeled $\la$) and $N\subset M$ be an $A$-submodule; we show that $N$ has a complement in $M$. Let $E:M\to N$ be a projection and define $E_0:M\to N$ by
\[ E_0(m)=\Lambda\o\la E(S\Lambda\t\la m)\]
for all $m\in M$. Then, for any $n\in N$,
\[ E_0(n)=\Lambda\o\la E(S\Lambda\t\la n)=\Lambda\o\la (S\Lambda\t\la n)=\eps(\Lambda)n=n.\]
It is easily seen that,
\[ \Delta(\Lambda)\tens a=a\o\o\Lambda\o\tens a\o\t\Lambda\t\tens a\t\]
by using the property $a\Lambda=\eps(a)\Lambda$,so for all $m\in M,a\in H$,
\begin{eqnarray*}
E_0(h\la m)	&	=	& \Lambda\o\la E(S\Lambda\t \la(a\la m))\\
	&	=	&	a\o\o\Lambda\o\la E(S(a\o\t\Lambda\t)\la(a\t\la m))\\
	&	=	&	a\o\o\Lambda\o\la E((S\Lambda\t)(Sa\o\t) a\t\la m)) \quad\text{since the action is associative}\\
	&	=	& a\Lambda\o\la E(S\Lambda\t\la m)\quad\text{by the Hopf coquasigroup axiom}\\
	&	=	& a\la E_0(m).
\end{eqnarray*}
So $E_0$ is a projection and an $A$-module map, therefore $\ker E_0$ is an $A$-complement for $N$ by a generalization of Maschke's theorem \cite{La69} and so $M$ is completely reducible. 
\endproof

\begin{remark}
To prove that a Hopf quasigroup $H$ is semisimple is not quite so straightforward. It seems that we may need extra conditions on $H$, or more specifically on the action of $H$ on a module $M$.
\end{remark}

\end{document}